\documentclass[12pt, eqno, twoside]{article}
\usepackage{graphicx}
\usepackage{amscd}
\usepackage[matrix,arrow,curve]{xy}
\usepackage{amssymb,amsmath}
\usepackage{amsmath,amsfonts,amssymb,latexsym}
\usepackage{fancyhdr}
\newtheorem{theo}{Theorem}[section]
\newtheorem{mydf}{Definition}[section]

\newtheorem{lem}{Lemma}[section]

\newenvironment{proof}[1][Proof]{\noindent\textbf{#1.} }{\
\rule{0.5em}{0.5em}}

\begin{document}

\pagestyle{fancy}
\fancyhead{} % clear all header fields
\fancyhead[EC]{\small\it Anastasios Mallios, \ Patrice P.
Ntumba}%
\fancyhead[EL,OR]{\thepage} \fancyhead[OC]{\small\it Pairings of
Sheaves of $\mathcal{A}$-Modules through Bilinear
$\mathcal{A}$-Morphisms }%
\fancyfoot{} % clear all footer fields
\renewcommand\headrulewidth{0.5pt}
\addtolength{\headheight}{2pt} % make space for the rule

\title{\Large{Pairings of Sheaves of $\mathcal{A}$-Modules through Bilinear
$\mathcal{A}$-Morphisms}}
\author{Anastasios Mallios, \ Patrice P.
Ntumba\footnote{Is the corresponding author for the paper}}
%wishes to thank the Department of Mathematics and Applied
%Mathematics, University of Pretoria, and the Algebra and Geometry
%Section, Department of Mathematics, University of Athens for
%sponsoring the research visit (01/12/2007--15/12/2007, University
%of Athens) during which this paper was finalized.}

\date{}
\maketitle

\begin{abstract}It is proved that for any free $\mathcal{A}$-modules $\mathcal{F}$ and $\mathcal{E}$ of finite rank on
some $\mathbb{C}$-algebraized space $(X, \mathcal{A})$ a
\textit{degenerate} bilinear $\mathcal{A}$-morphism $\Phi:
\mathcal{F}\times \mathcal{E}\longrightarrow \mathcal{A}$  induces a
\textit{non-degenerate} bilinear $\mathcal{A}$-morphism
$\overline{\Phi}: \mathcal{F}/\mathcal{E}^\perp\times
\mathcal{E}/\mathcal{F}^\perp\longrightarrow \mathcal{A}$, where
$\mathcal{E}^\perp$ and $\mathcal{F}^\perp$ are the
\textit{orthogonal} sub-$\mathcal{A}$-modules associated with
$\mathcal{E}$ and $\mathcal{F}$, respectively. This result
generalizes the finite case of the classical result, which states
that given two vector spaces $W$ and $V$, paired into a field $k$,
the induced vector spaces $W/V^\perp$ and $V/W^\perp$ have the same
dimension. Some related results are discussed as well.
\end{abstract}

\maketitle

\section*{Introduction} The goal of this paper is to provide
additional steps for developing \textit{classical symplectic
geometry} within the setting of \textit{Abstract Differential
Geometry} (\textit{ADG} in short) ($\grave{a}$ la Mallios), cf.
\cite{mallios} and \cite{mallios1}. The attempt of taking ADG to new
horizons, such those related to classical symplectic geometry was
initiated in the our paper \cite{malliosntumba}. The main result in
the article \cite{malliosntumba} is that given an ordered
$\mathbb{R}$-algebraized space $(X, \mathcal{A}, \mathcal{P}, |\cdot
|)$, the standard free $\mathcal{A}$-module $\mathcal{A}^n$ of rank
$n$ on $X$, there exists for every open subset $U$ of $X$ a basis
$\mathcal{B}(U)$ of $\mathcal{A}^n(U)$, relative to which the matrix
of a non-zero skew-symmetric and non-degenerate bilinear sheaf
$\mathcal{A}$-morphism $\omega\equiv (\omega_U): \mathcal{A}^n\oplus
\mathcal{A}^n\longrightarrow \mathcal{A}$ is the matrix
\[\left(\begin{array}{cc} 0 & \mbox{I}_n\\ -\mbox{I}_n &
0\end{array}\right).\]

In order to aptly pursue the goal of our ongoing research project,
we are laying some ground work, regarding pairings of sheaves of
$\mathcal{A}$-modules. Duality and pairings of
$\mathcal{A}$-modules, as we will see in our subsequent work, are a
necessary prerequisite for the layout of Abstract Symplectic
Geometry. In this article, the most important results are contained
in Theorems \ref{theo2} and \ref{theo3}. Theorem \ref{theo2}
contends that given $\mathcal{A}$-modules $\mathcal{F}$ and
$\mathcal{E}$, paired into $\mathcal{A}$ such that the left kernel
$\mathcal{E}^\perp=0$, one may find natural
$\mathcal{A}$-isomorphisms into \[\begin{array}{ll}
\mathcal{E}/\mathcal{F}_0^\perp\longrightarrow \mathcal{F}_0^\ast, &
\mbox{and $ \mathcal{E}_0^\perp\longrightarrow
(\mathcal{E}/\mathcal{E}_0)^\ast$}\end{array}\]for
sub-$\mathcal{A}$-modules $\mathcal{F}_0$ and $\mathcal{E}_0$ of
$\mathcal{F}$ and $\mathcal{E}$, respectively. In the special case
where $\mathcal{F}$ and $\mathcal{E}$ are free $\mathcal{A}$-modules
of finite rank, if $\Phi: \mathcal{F}\times
\mathcal{E}\longrightarrow \mathcal{A}$ is bilinear, then
$\mathcal{F}$ is $\mathcal{A}$-isomorphic to $\mathcal{E}$ provided
that $\Phi$ is non-degenerate, and if $\Phi$ is degenerate, it
induces a non-degenerate $\mathcal{A}$-morphism $$\overline{\Phi}:
\mathcal{F}/\mathcal{E}^\perp\times
\mathcal{E}/\mathcal{F}^\perp\longrightarrow \mathcal{A}$$ such that
\[(\mathcal{F}/\mathcal{E}^\perp)^\perp=0 =
(\mathcal{E}/\mathcal{F}^\perp)^\perp.\]

\section{Theorems on ranks of free $\mathcal{A}$-modules}Let $\mathcal{E}\equiv (\mathcal{E},
\pi, X)$ and $\mathcal{F}\equiv (\mathcal{F}, \rho, X)$ be two
$\mathcal{A}$-modules on a topological space $X$, and let
\[\Gamma(\mathcal{E})\equiv (\Gamma(\mathcal{E})(U)\equiv \Gamma(U,
\mathcal{E})\equiv \mathcal{E}(U), \pi^U_V)\]and
\[\Gamma(\mathcal{F})\equiv (\mathcal{F}(U), \rho^U_V)\] be the
complete presheaves of sections of $\mathcal{E}$ and
$\mathcal{F}$, respectively. For every open set $U\subseteq X$,
\[(\mathcal{E}+ \mathcal{F})(U)\equiv \mathcal{E}(U)+
\mathcal{F}(U):= \Gamma(U, \mathcal{E})+ \Gamma(U,
\mathcal{F})\]is the sum of $\mathcal{A}(U)$-modules
$\mathcal{E}(U)$ and $\mathcal{F}(U)$, where
$\Gamma(\mathcal{A})\equiv (\Gamma(U, \mathcal{A}), \kappa^U_V)$
is the complete presheaf of sections of the coefficient sheaf
$\mathcal{A}$. Next, consider the following presheaf on $X$,
defined by the rule \begin{equation}\label{eq1} U\longmapsto
S(U):= \mathcal{E}(U)+ \mathcal{F}(U),\end{equation}where $U$ is
open in $X$ and restriction maps are maps $\lambda^U_V:
S(U)\longmapsto S(V)$ such that \[\lambda^U_V(s+ t)= \pi^U_V(s)+
\rho^U_V(t),\]for all $s\in \mathcal{E}(U)$ and $t\in
\mathcal{F}(U)$. Thus, (\ref{eq1}) yields a
$\Gamma(\mathcal{A})$-presheaf on $X$, which we denote \[S\equiv
\Gamma(\mathcal{E}+ \mathcal{F}):= \Gamma(\mathcal{E})+
\Gamma(\mathcal{F}).\]The $\mathcal{A}$-module \[\mathcal{E}+
\mathcal{F}:= \mathbf{S}(S)\equiv \mathbf{S}(\Gamma(\mathcal{E})+
\Gamma(\mathcal{F}))\]generated by the presheaf $S$ is called the
\textit{sum} of $\mathcal{A}$-modules $\mathcal{E}$ and
$\mathcal{F}$.

The reader will have no difficulty in proving that the sheaf
$\mathcal{E}+ \mathcal{F}$ is indeed an $\mathcal{A}$-module, and
$S$ is complete.

\begin{theo}\label{theo1}
Let $\mathcal{E}$ be an $\mathcal{A}$-module on a topological
space $X$, $\mathcal{F}$ and $\mathcal{G}$ be
sub-$\mathcal{A}$-modules of $\mathcal{E}$. Then, \[(\mathcal{F}+
\mathcal{G})/\mathcal{F}= \mathcal{G}/(\mathcal{F}\cap
\mathcal{G})\]within an $\mathcal{A}$-isomorphism.
\end{theo}

\begin{proof}
Let $\phi\equiv (\phi_U)_{U\in \tau}$ be the (canonical) quotient
$\Gamma( \mathcal{A})$-morphism $$\Gamma(\mathcal{E})\longrightarrow
\Gamma(\mathcal{E})/\Gamma(\mathcal{F}),$$ where
$\Gamma(\mathcal{E})\equiv (\mathcal{E}(U), \pi^U_V)$ is the
complete presheaf of sections of $\mathcal{E}$, and
\[\Gamma(\mathcal{E})/\Gamma(\mathcal{F})\equiv
((\Gamma(\mathcal{E})/\Gamma(\mathcal{F}))(U)\equiv \Gamma(U,
\mathcal{E})/\Gamma(U, \mathcal{F})\equiv
\mathcal{E}(U)/\mathcal{F}(U), \sigma^U_V)\](the $\sigma^U_V$ are
the obvious restriction maps, given by \[\sigma^U_V(s+ \Gamma(U,
\mathcal{F}))= \pi^U_V(s)+ \Gamma(V, \mathcal{F}),\]for all $s\in
\Gamma(U, \mathcal{E})$) is the generating presheaf of the quotient
$\mathcal{A}$-module $\mathcal{E}/\mathcal{F}$, cf.
Mallios[\cite{mallios}, pp. 114, 115]. For every open $U\subseteq
X$, the restriction $\psi_U$ of $\phi_U$ to the
sub-$\mathcal{A}(U)$-module $\Gamma(U, \mathcal{G})$ is the
canonically constructed $\mathcal{A}(U)$-morphism \[\psi_U:
\mathcal{G}(U)\longrightarrow \mathcal{E}(U)/\mathcal{F}(U),\]given
by
\[\begin{array}{ll}\psi_U(s):= s+ \mathcal{F}(U), & s\in
\mathcal{G}(U).\end{array}\]For every fixed open $U\subseteq X$,
the union of the cosets $s+ \mathcal{F}(U)$, $s\in
\mathcal{G}(U)$, forms the $\mathcal{A}(U)$-module
$\mathcal{G}(U)+ \mathcal{F}(U)\equiv (\mathcal{G}+
\mathcal{F})(U)$; therefore \[\psi_U(\mathcal{G}(U))=
(\mathcal{G}(U)+ \mathcal{F}(U))/\mathcal{F}(U)\equiv
(\mathcal{G}+ \mathcal{F})(U)/\mathcal{F}(U).\]But, for all
elements $s\in \mathcal{G}(U)$, with $U$ as above, i.e. an open
set in $X$, we have \[\psi_U(s)= \phi_U(s)\]and \[\ker \phi_U=
\mathcal{F}(U).\]Therefore, \[\ker \psi_U= \mathcal{F}(U)\cap
\mathcal{G}(U)\equiv (\mathcal{F}\cap \mathcal{G})(U).\]By
elementary algebra, we construct the canonical
$\mathcal{A}(U)$-isomorphism \[\overline{\psi}_U:
\mathcal{G}(U)/(\mathcal{F}\cap \mathcal{G})(U)\longrightarrow
(\mathcal{F}+ \mathcal{G})(U)/\mathcal{F}(U),\]for every open
$U\subseteq X$. If $V\subseteq U$ is an open subset of a certain
fixed open $U\subseteq X$, it is not difficult to see that the
diagram \[\xymatrix{\mathcal{G}(U)/(\mathcal{F}\cap
\mathcal{G})(U)\ar[r]^{\overline{\psi}_U}\ar[d] & (\mathcal{F}+
\mathcal{G})(U)/\mathcal{F}(U)\ar[d] \\
\mathcal{G}(V)/(\mathcal{F}\cap
\mathcal{G})(V)\ar[r]^{\overline{\psi}_V} & (\mathcal{F}+
\mathcal{G})(V)/\mathcal{F}(V) }\]commutes. Consequently,
$\overline{\psi}\equiv (\overline{\psi}_U)_{U\in \tau}$ is a
$\Gamma(\mathcal{A})$-isomorphism of the presheaves
$\Gamma(\mathcal{G})/\Gamma(\mathcal{F}\cap \mathcal{G})$ and
$\Gamma(\mathcal{F}+ \mathcal{G})/\Gamma(\mathcal{F})$. Applying
the sheafification functor $\mathbf{S}$ (see
Mallios[\cite{mallios}, p. 33]) to the diagram
\[\overline{\psi}\equiv (\overline{\psi}_U): \Gamma(\mathcal{G})/\Gamma(\mathcal{F}\cap
\mathcal{G})\longrightarrow \Gamma(\mathcal{F}+
\mathcal{G})/\Gamma(\mathcal{F}),\] we obtain the sought
$\mathcal{A}$-isomorphism \[\overline{\psi}:
\mathcal{G}/(\mathcal{F}\cap \mathcal{G})\longrightarrow
(\mathcal{F}+ \mathcal{G})/\mathcal{F}.\]
\end{proof}

In the special case that $\mathcal{E}= \mathcal{F}\oplus
\mathcal{G}$, it follows for any open $U\subseteq X$, cf.
Mallios[\cite{mallios}, pp. 121, 122], that \[\Gamma(U,
\mathcal{E})= \Gamma(U, \mathcal{F})\oplus \Gamma(U,
\mathcal{G}).\]Applying Artin[\cite{artin}, Theorem 1.2, p. 7], we
obtain that $\Gamma(U, \mathcal{E})/\Gamma(U, \mathcal{F})$ is
$\mathcal{A}(U)$-isomorphic to $\Gamma(U, \mathcal{G})/0=
\Gamma(U, \mathcal{G}$, where $U$ is any arbitrary open set in
$X$. It thus follows that \[\mathcal{E}/\mathcal{F}:=
\textbf{S}(\Gamma(\mathcal{E})/\Gamma(\mathcal{F}))=
\textbf{S}(\Gamma(\mathcal{G}))= \mathcal{G},\]that is
\[\mathcal{E}/\mathcal{F}= \mathcal{G}\]within an
$\mathcal{A}$-isomorphism.

Before we proceed to some more theorems on $\mathcal{A}$-modules,
let us define what is meant by \textit{free
sub-$\mathcal{A}$-modules of a free $\mathcal{A}$-module}.

\begin{mydf}
\emph{Let $\mathcal{E}$ be the free $\mathcal{A}$-module
$\mathcal{A}^{(I)}:= \oplus_I\mathcal{A}$, where $I$ is an arbitrary
indexing set, and let $\mathcal{F}\subseteq \mathcal{E}$ be a
sub-$\mathcal{A}$-module of $\mathcal{E}$ such that
\[\mathcal{F}_x:= \mathcal{A}^{(J)}_x\oplus
\underbrace{0\oplus \ldots\oplus 0}_{I\setminus J}\subseteq
\mathcal{A}^{(I)}_x=: \mathcal{E}_x\] for all $x\in X$, and where
$J$ is a subset of $I$. $\mathcal{F}$ is called a \textit{free
sub-$\mathcal{A}$-module} of $\mathcal{E}$, and is easily identified
with the free $\mathcal{A}$-module $\mathcal{A}^{(J)}$. The free
sub-$\mathcal{A}$-module $\mathcal{G}:= \mathcal{A}^{(I\setminus
J)}$ is called a free sub-$\mathcal{A}$-module \textit{supplementary
to} $\mathcal{F}$. It is obvious that $\mathcal{E}=
\mathcal{F}\oplus \mathcal{G}$. The fibers of $\mathcal{G}$ are
$\mathcal{A}_x$-modules \[\mathcal{G}_x=
\underbrace{0\oplus\ldots\oplus 0}_I\oplus \mathcal{A}^{(I\setminus
J)}_x\subseteq \mathcal{A}^I_x. \]} \hfill$\square$
\end{mydf}

Theorem 1.4, cf. Artin[\cite{artin}, p. 9], is also immediate. In
effect, let $\mathcal{E}$ be a free $\mathcal{A}$-module on a
topological space $X$, viz. $\mathcal{E}\cong \mathcal{A}^{(I)}$,
where $I$ is arbitrary, and let $\mathcal{F}$ be a free
sub-$\mathcal{A}$-module of $\mathcal{E}$. If $\mathcal{G}$ is a
free sub-$\mathcal{A}$-module supplementary to $\mathcal{F}$ in
$\mathcal{E}$, we have that \[\mbox{rank}\ \mathcal{E}= \mbox{rank}
\ \mathcal{F}+ \mbox{rank}\ \mathcal{G},\]and since
$\mathcal{G}\cong \mathcal{E}/\mathcal{F}$, one obtains
\[\mbox{rank}\ \mathcal{E}= \mbox{rank}\ \mathcal{F}+ \mbox{rank}\
\mathcal{E}/\mathcal{F}.\]

Let now $\mathcal{F}_1\subseteq \mathcal{F}_2\subseteq
\mathcal{F}_3$ be free sub-$\mathcal{A}$-modules of a free
$\mathcal{A}$-module $\mathcal{E}$. For some free
sub-$\mathcal{A}$-modules $\mathcal{G}_2$ and $\mathcal{G}_3$ of
$\mathcal{E}$, we have \[\begin{array}{ll} \mathcal{F}_2=
\mathcal{F}_1\oplus \mathcal{G}_2, & \mathcal{F}_3=
\mathcal{F}_2\oplus \mathcal{G}_3\end{array}\]and, therefore,
\[\mathcal{F}_3= \mathcal{F}_1\oplus (\mathcal{G}_2\oplus
\mathcal{G}_3).\]

By virtue of Theorem \ref{theo1}, $\mathcal{G}_2$, $\mathcal{G}_3$
and $\mathcal{G}_2\oplus \mathcal{G}_3$ are
$\mathcal{A}$-isomorphic to $\mathcal{F}_2/\mathcal{F}_1$,
$\mathcal{F}_3/ \mathcal{F}_2$ and $\mathcal{F}_3/\mathcal{F}_1$,
respectively; consequently \begin{eqnarray*} \mbox{rank}\
\mathcal{G}_2 & = & \mbox{rank}\ \mathcal{F}_2/\mathcal{F}_1\\
\mbox{rank}\ \mathcal{G}_3 & = & \mbox{rank}\
\mathcal{F}_3/\mathcal{F}_2\\ \mbox{rank}\ (\mathcal{G}_2\oplus
\mathcal{G}_3) & = & \mbox{rank}\
\mathcal{F}_3/\mathcal{F}_1.\end{eqnarray*}Thus, we obtain
\begin{equation} \label{eq2}\mbox{rank}\
\mathcal{F}_3/\mathcal{F}_1= \mbox{rank}\
\mathcal{F}_2/\mathcal{F}_1+ \mbox{rank}\
\mathcal{F}_3/\mathcal{F}_2.\end{equation}

Let now $\mathcal{F}$ and $\mathcal{G}$ be two given free
sub-$\mathcal{A}$-modules of a free $\mathcal{A}$-module
$\mathcal{E}$. For $\mathcal{F}_1= 0$, $\mathcal{F}_2=
\mathcal{F}$, $\mathcal{F}_3= \mathcal{F}+ \mathcal{G}$, Equation
(\ref{eq2}) becomes \begin{eqnarray*} \mbox{rank}\ (\mathcal{F}+
\mathcal{G})& = & \mbox{rank}\ \mathcal{F}+ \mbox{rank}\
(\mathcal{F}+ \mathcal{G})/\mathcal{F}\\ & = & \mbox{rank}\
\mathcal{F}+ \mbox{rank}\ \mathcal{G}/(\mathcal{F}\cap
\mathcal{G}). \end{eqnarray*} Adding $\mbox{rank}\
(\mathcal{F}\cap \mathcal{G})$ to both sides of the last equation,
and using the fact that \[\mbox{rank}\ \mathcal{G}= \mbox{rank}\
(\mathcal{F}\cap \mathcal{G})+ \mbox{rank}\
\mathcal{G}/(\mathcal{F}\cap \mathcal{G}),\]we obtain
\[\mbox{rank}\ (\mathcal{F}+ \mathcal{G})+ \mbox{rank}\
(\mathcal{F}\cap \mathcal{G})= \mbox{rank}\ \mathcal{F}+
\mbox{rank}\ \mathcal{G}.\] Next, we put in Equation (\ref{eq2})
$\mathcal{F}_1= \mathcal{F}\cap \mathcal{G}$, $\mathcal{F}_2=
\mathcal{G}$, $\mathcal{F}_3= \mathcal{E}$ to get
\begin{eqnarray*} \mbox{rank}\ \mathcal{E}/(\mathcal{F}\cap
\mathcal{G}) & = & \mbox{rank}\ \mathcal{G}/(\mathcal{F}\cap
\mathcal{G})+ \mbox{rank}\ \mathcal{E}/\mathcal{G}\\ & = &
\mbox{rank}\ (\mathcal{F}+ \mathcal{G})/\mathcal{F}+ \mbox{rank}\
\mathcal{E}/\mathcal{G}.\end{eqnarray*}

If we add $\mbox{rank}\ \mathcal{E}/(\mathcal{F}+ \mathcal{G})$,
and use the equation \[\mbox{rank}\ \mathcal{E}/\mathcal{F}=
\mbox{rank}\ (\mathcal{F}+ \mathcal{G})/\mathcal{F}+ \mbox{rank}\
\mathcal{E}/(\mathcal{F}+ \mathcal{G})\](Put $\mathcal{F}_1=
\mathcal{F}$, $\mathcal{F}_2= \mathcal{F}+ \mathcal{G}$,
$\mathcal{F}_3= \mathcal{E}$ in Equation (\ref{eq2}) to obtain the
last equation.), we obtain
\[\mbox{rank}\ \mathcal{E}/(\mathcal{F}+
\mathcal{G})+ \mbox{rank}\ \mathcal{E}/(\mathcal{F}\cap
\mathcal{G}) =\mbox{rank}\ \mathcal{E}/\mathcal{F}+ \mbox{rank}\
\mathcal{E}/\mathcal{G}.\]

\begin{mydf}
\emph{Let $\mathcal{E}$ be a free $\mathcal{A}$-module on a
topological space $X$, and $\mathcal{F}$ a free
sub-$\mathcal{A}$-module of $\mathcal{E}$ such that its supplement
is a free sub-$\mathcal{A}$-module of \textit{finite rank}. The
rank of the free sub-$\mathcal{A}$-module
$\mathcal{E}/\mathcal{F}$ is called the \textit{corank} of
$\mathcal{F}$, viz. \[\mbox{corank}\ \mathcal{F}= \mbox{rank}\
\mathcal{E}/\mathcal{F}.\]} \hfill$\square$
\end{mydf}

The above various results can be expressed as follows:

\begin{theo}Let $\mathcal{E}$ be a free $\mathcal{A}$-module on a
topological space $X$, and $\mathcal{F}$ and $\mathcal{G}$ free
sub-$\mathcal{A}$-modules of $\mathcal{E}$. Then,
\begin{eqnarray*} \mbox{rank}\ \mathcal{F}+ \mbox{corank}\
\mathcal{F} & = & \mbox{rank}\ \mathcal{E}\\ \mbox{rank}\
(\mathcal{F}+ \mathcal{G})+ \mbox{rank}\ (\mathcal{F}\cap
\mathcal{G}) & = & \mbox{rank}\ \mathcal{F}+ \mbox{rank}\
\mathcal{G}\\ \mbox{corank}\ (\mathcal{F}+ \mathcal{G})+
\mbox{corank}\ (\mathcal{F}\cap \mathcal{G}) & = & \mbox{corank}\
\mathcal{F}+ \mbox{corank}\ \mathcal{G}.\end{eqnarray*}
\end{theo}

\section{Pairings}

\begin{mydf}
\emph{Let $(X, \mathcal{A})$ be a $\mathbb{C}$-algebraized space,
and let $\mathcal{E}$ and $\mathcal{F}$ be $\mathcal{A}$-modules
on $X$. We say that $\mathcal{F}$ and $\mathcal{E}$ are
\textbf{paired into} $\mathcal{A}$, or $[\mathcal{F}, \mathcal{E};
\mathcal{A}]$ is a \textbf{pairing} provided a \textit{bilinear
$\Gamma(\mathcal{A})$-morphism} $\vartheta:
\Gamma(\mathcal{F})\oplus \Gamma(\mathcal{E})\longrightarrow
\Gamma(\mathcal{A})$ is defined between the
$\Gamma(\mathcal{A})$-presheaves $\Gamma(\mathcal{F})\oplus
\Gamma(\mathcal{E})$ and $\Gamma(\mathcal{A})$.\hfill$\square$}
\end{mydf}

We notice for any $\mathcal{A}$-module $\mathcal{E}$ on $X$,
$\mathcal{E}^\ast$ and $\mathcal{E}$ are paired into $\mathcal{A}$
by the bilinear $\Gamma{\mathcal{A}}$-morphism $\vartheta:
\Gamma(\mathcal{E}^\ast)\oplus \Gamma(\mathcal{E})\longrightarrow
\Gamma(\mathcal{A})$, given by
\[\vartheta_U(\phi, s):= \phi(s)\]for $\phi\equiv (\phi_V)_{U\supseteq
V,\ open}\in \mathcal{E}^\ast(U):=
\mathcal{H}om_\mathcal{A}(\mathcal{E}, \mathcal{A})(U)\equiv
Hom_{\mathcal{A}|_U}(\mathcal{E}|_U, \mathcal{A}|_U)$ (cf.
Mallios[\cite{mallios}, relation (5.1), p. 298, Definition 6.1,
p.134]), $s\in \mathcal{E}(U)$, and $U$ an open subset of $X$.

Now, consider a free $\mathcal{A}$-module $\mathcal{E}$, i.e.
$\mathcal{E}= \mathcal{A}^{(I)}$ within an
$\mathcal{A}$-isomorphic, and let $\{s_i\}$ be a basis of the
$\mathcal{A}(U)$-module $\mathcal{E}(U)$, where $U$ is an open
subset of $X$. (Cf. Mallios[\cite{mallios}, (3.13), p.122] For any
open $U\subseteq X$, $\mathcal{E}(U)$ is
$\mathcal{A}(U)$-isomorphic to $\prod_{i\in I}\mathcal{A}(U)\equiv
\prod_{i\in I}\Gamma(U, \mathcal{A})$.) Let $\phi\equiv
(\phi_V)_{U\supseteq V,\ open}$ be an element of
$\mathcal{E}^\ast(U)$; put
\[\vartheta(\phi, s_i)\equiv \phi_U(s_i):= a_i\in \mathcal{A}(U),\]for all
$i\in I$. An element $s\in \mathcal{E}(U)$ is written uniquely in
the form $s= \sum_{i\in I}s_ir_i$, where $0\neq r_i\in
\mathcal{A}(U)$ holds only for finitely many indices $i$. Then,
\[\vartheta(\phi, s)\equiv \phi_U(s)= \sum_{i\in I}\phi_U(s_i)r_i=: \sum_{i\in
I}a_ir_i.\]

Now, assume that the restriction maps for the presheaf of sections
of $\mathcal{E}$ (resp. $\mathcal{A}$) are given by $\rho^U_V$
(resp. $\sigma^U_V$), where $V$ and $U$ are open subsets of $X$,
with $V\subseteq U$. For all $s\in \mathcal{E}(U)$, we have for
any open $V\subseteq U$,
\begin{eqnarray*}\phi_V(\rho^U_V(s))\equiv \phi_V(s|_V)=
\sigma^U_V(\phi_U(s))= \sigma^U_V(\sum_{i\in I}a_ir_i)=:
\sum_{i\in I}a_i|_{V}r_i|_V;\end{eqnarray*}moreover, since
$\{\rho^U_V(s_i)\}_{i\in I}$ is a basis of $\mathcal{E}(V)$, it
follows that $\phi\equiv (\phi_V)\in \mathcal{E}(U)$ is known if
all the $a_i\in \mathcal{A}(U)$ are known.

Select conversely an $a_i\in \mathcal{A}(U)$ for each index $i$,
and define \[\vartheta_U(\phi, s):= \phi_U(s)= \sum_{i\in
I}a_ir_i,\]assuming as above that $\phi\equiv (\phi_V)_{U\supseteq
V,\ open}\in \mathcal{E}^\ast(U)$, $s\in \mathcal{E}(U)$ and $s=
\sum_{i\in I}s_ir_i$, where $\{s_i\}_{i\in I}$ is a basis of
$\mathcal{E}(U)$. The sum $\sum_{i\in I}a_ir_i$ is finite since
only finitely many $r_i$ are non-zero. That $\phi_U(s+t)=
\phi_U(s)+ \phi_U(t)$, and $\phi_U(sa)= \phi_U(s)a$ for all $s,
t\in \mathcal{E}(U)$, and $a\in \mathcal{A}(U)$ is immediately
clear. Since $s_i= \sum_{j\in I}s_j\delta_{ji}$, where, as usual
$\delta_{ii}= 1\in \mathcal{A}(U)$ and $\delta_{ji}= 0$ for $j\neq
i$, we obtain $\vartheta_U(\phi, s_i):= \phi_U(s_i)= \sum_{j\in
I}a_j\delta_{ji}= a_i.$ Thus,

\begin{theo}
Let $(X, \mathcal{A})$ be a $\mathbb{C}$-algebraized space and
$\mathcal{E}$ be a free $\mathcal{A}$-module. If $\{s_i\}$ is a
basis of $\mathcal{E}(U)$, where $U$ is open in $X$, then for
arbitrarily chosen sections $a_i\in \mathcal{A}(U)$, \textsf{there
is one and only one} $\phi\in \mathcal{E}^\ast(U)$ such that
\[\phi_U(s_i)= a_i.\]
\end{theo}

As above let $\mathcal{E}$ be a free $\mathcal{A}$-module on a
topological space $X$, viz. $\mathcal{E}= \mathcal{A}^{(I)}$ within
an $\mathcal{A}$-isomorphism, and let $\{s_i\}$ be a basis of
$\mathcal{E}(U)= \mathcal{A}^{(I)}(U):= \Gamma(U,
\mathcal{A}^{(I)})= \Gamma(U, \mathcal{A})^{(I)}.$ Denote by
$\phi_i$ the $\mathcal{A}|_U$-morphism in $\mathcal{E}^\ast(U):=
Hom_{\mathcal{A}|_U}(\mathcal{E}|_U, \mathcal{A}|_U)$ for which
$\phi_i(s_j):= \phi_{i,U}(s_j)= \delta_{ij, U}.$ Let $B(U)$ be the
sub-$\mathcal{A}|_U$-module of $\mathcal{E}^\ast(U)$, spanned by the
$\mathcal{A}|_U$-morphisms $\phi_i$. That is, given $\phi\in
\mathcal{E}^\ast(U)$, $\phi\in B(U)$ provided $\phi= \sum_{i\in
I}\alpha_i\phi_i$, where $\alpha_i\in \mathcal{A}|_U$ for all $i\in
I$, and only finitely many $\alpha_i$ are non-zero. Since there are
as many $\phi_i$ as there are $s_i$, we get that $\dim B(U)= \dim
\mathcal{E}(U)$. If $I$ is finite, $\mathcal{E}^\ast= \mathcal{E}$
within an $\mathcal{A}$-isomorphism (cf. Mallios[\cite{mallios}, p.
298]); therefore $\dim \mathcal{E}^\ast(U)= \dim \mathcal{E}(U)$ for
any open subset $U\subseteq X$. If $\dim \mathcal{E}(U)= \infty$,
then $\dim B(U)= \infty$ and, since $B(U)\subseteq
\mathcal{E}^\ast(U)$, we may put $\dim \mathcal{E}^\ast(U)= \infty$.
So in this case as well, we have \[\dim \mathcal{E}(U)= \dim
\mathcal{E}^\ast(U),\]for any open $U\subseteq X$. Next, if $s=
\sum_{i\in I}s_ia_i\in \mathcal{E}(U)$, where $U$ is a fixed open
subset of $X$, and $s\neq 0$, then at least one $a_j\neq 0$, so that
$\phi_j(s)\neq 0$. By definition of $\mathcal{E}^\ast$, we know
trivially that only the zero section vanishes on all of
$\mathcal{E}(U)$. Now, we see an analogue: If $s\in \mathcal{E}(U)$,
and $\phi(s)= 0$ for all $\phi\in \mathcal{E}^\ast(U)$, then $s=0$.

Hence, we have

\begin{theo}
Let $\mathcal{E}$ be a free $\mathcal{A}$-module on a topological
space $X$. Then, for any open subset $U\subseteq X$, $\dim
\mathcal{E}^\ast(U)= \dim \mathcal{E}(U)$. If $\phi(s)= 0$ for all
$s\in \mathcal{E}(U)$, then $\phi= 0$; on the other side, if
$\phi(s)=0$ for all $\phi\in \mathcal{E}^\ast(U)$, then $s=0$.
Finally, let $\dim \mathcal{E}(U)= n$ for some fixed open
$U\subseteq X$, then $\dim \mathcal{E}(V)= n$ for any open
$V\subseteq X$. To a given basis $\{s_i\}$ of $\mathcal{E}(U)$, we
can find a \textsf{dual basis} $\{\phi_i\}$ of $\mathcal{E}^\ast(U)=
\mathcal{E}(U)$, where \[\phi_i(s_j):= \phi_{i, U}(s_j)= \delta_{ij,
U}\in \mathcal{A}(U).\]
\end{theo}

Turning over to pairings of $\mathcal{A}$-modules, we suppose that
$\mathcal{A}$-modules $\mathcal{F}$ and $\mathcal{E}$, defined on a
topological space $X$, are given and form a pairing into
$\mathcal{A}$.

\begin{mydf}
\emph{Let $U$ be an open subset of $X$, $t\in \mathcal{F}(U)$ and
$s\in \mathcal{E}(U)$. We say that $t$ is \textbf{orthogonal} to $s$
provided if $\vartheta\equiv (\vartheta_U)_{X\supseteq U,\ open}:
\Gamma(\mathcal{F})\oplus \Gamma(\mathcal{E})\longrightarrow \Gamma
\mathcal{A}$ the bilinear $\Gamma(\mathcal{A})$-morphism defining
the pairing $[\mathcal{F}, \mathcal{E}: \mathcal{A}]$, then
$\vartheta_U(t, s):= ts= 0$, i.e. $t(x)s(x)= 0$, for all $x\in U$.
More generally, $\mathcal{F}(U)$ is said to be \textit{orthogonal}
to $\mathcal{E}(U)$ if $ts= 0$ for all $t\in \mathcal{F}(U)$ and
$s\in \mathcal{E}(U)$. Similarly, a sub-$\mathcal{A}$-module
$\mathcal{F}_0$ of $\mathcal{F}$ is \textit{orthogonal} to a
sub-$\mathcal{A}$-module $\mathcal{E}_0$ of $\mathcal{E}$ if
$\mathcal{F}_0(U)$ is orthogonal to $\mathcal{E}_0(U)$ for any open
set $U\subseteq X$.}\hfill$\square$
\end{mydf}

\begin{lem}
Let $[\mathcal{F}, \mathcal{E}; \mathcal{A}]$ be a pairing in which
$\mathcal{F}$ and $\mathcal{E}$ are $\mathcal{A}$-modules on a
topological space $X$, and let $\mathcal{F}_0$ and $\mathcal{E}_0$
be sub-$\mathcal{A}$-modules of $\mathcal{F}$ and $\mathcal{E}$,
respectively. Furthermore, let \[\mathcal{F}_0(U)^\perp =\{s\in
\mathcal{E}(U):\ ts=0\ \mbox{for all $t\in \mathcal{F}_0(U)\subseteq
\mathcal{F}(U)$}\}\]and \[\mathcal{E}_0(U)^\perp = \{t\in
\mathcal{F}(U):\ ts=0\ \mbox{for all $s\in \mathcal{E}_0(U)\subseteq
\mathcal{E}(U)$}\}\]for all open $U\subseteq X$, and let
$(\mathcal{F}(U), \rho^U_V)$ and $(\mathcal{E}(U), \pi^U_V)$ be
$($complete$)$ presheaves of sections of $\mathcal{F}$ and
$\mathcal{E}$, respectively. The \textsf{sheaf} generated by the
presheaf, given by the correspondence \[\begin{array}{ll}
U\longmapsto \mathcal{E}_0(U)^\perp & (\mbox{resp. $U\longmapsto
\mathcal{F}_0(U)^\perp$}),\end{array}\]where $U$ is an open subset
of $X$, along with restriction maps
\[\begin{array}{ll}(\pi_0^\perp)^U_V:
\mathcal{E}_0(U)^\perp\longrightarrow \mathcal{E}_0(V)^\perp &
(\mbox{resp. $(\rho_0^\perp)^U_V:
\mathcal{F}_0(U)^\perp\longrightarrow
\mathcal{F}_0(V)^\perp$})\end{array}\] such that \[\begin{array}{ll}
(\pi_0^\perp)^U_V:= \rho^U_V|_{\mathcal{E}_0(U)^\perp} &
(\mbox{resp. $(\rho_0^\perp)^U_V:=
\pi^U_V|_{\mathcal{F}_0(U)^\perp}$})\end{array}\] is a
\textsf{sub-$\mathcal{A}$-module} of $\mathcal{F}$ $($resp.
$\mathcal{E}$ $)$, and is called the
\textsf{sub-$\mathcal{A}$-module orthogonal to} $\mathcal{E}_0$
$($resp. $\mathcal{F}_0$ $)$. We will denote by \[\begin{array}{ll}
\mathcal{E}_0^\perp & (\mbox{resp.
$\mathcal{F}_0^\perp$})\end{array}\]the sub-$\mathcal{A}$-module
orthogonal to $\mathcal{E}_0$ $($resp. $\mathcal{F}_0$ $)$, thus
obtained.
\end{lem}

\begin{proof}
For any open $U\subseteq X$, one sees easily that
$\mathcal{E}_0(U)^\perp$ and $\mathcal{F}_0(U)^\perp$ are
sub-$\mathcal{A}(U)$-modules of $\mathcal{A}(U)$-modules
$\mathcal{F}(U)$ and $\mathcal{E}(U)$, respectively. It follows from
Mallios-Ntumba[\cite{malliosntumba1}, Definition 1.1] that
$\mathcal{E}_0^\perp$ and $\mathcal{F}_0^\perp$ are subsheaves of
$\mathcal{F}$ and $\mathcal{E}$, respectively. Finally, the
sheafifications $\mathcal{E}_0^\perp$ and $\mathcal{F}_0^\perp$, by
virtue of Mallios[\cite{mallios}, Statement (1.54), p. 101], are
sub-$\mathcal{A}$-modules of $\mathcal{F}$ and $\mathcal{E}$,
respectively.
\end{proof}

Lemma 1.2 in Mallios-Ntumba\cite{malliosntumba1} is a particular
case of Lemma \ref{lem1.2} below; the proof of Lemma 1.2
\cite{malliosntumba1} applies here as well.

\begin{lem}\label{lem1.2}
Let $\mathcal{E}$ and $\mathcal{F}$ be $\mathcal{A}$-modules on a
topological space $X$, and suppose that $[\mathcal{F}, \mathcal{E};
\mathcal{A}]$ is a pairing. Then, for any sub-$\mathcal{A}$-module
$\mathcal{E}_0$ of $\mathcal{E}$, the correspondence \[U\longmapsto
\mathcal{E}_0(U)^\perp\]along with maps $(\pi_0^\perp)^U_V$, as
defined above, yields a complete presheaf of $\mathcal{A}$-modules
on $X$. Similarly, for any sub-$\mathcal{A}$-module $\mathcal{F}_0$
of $\mathcal{F}$, the assignment \[U\longmapsto
\mathcal{F}_0(U)^\perp\]with the afore-defined maps
$(\rho_0^\perp)^U_V$ define a complete presheaf of
$\mathcal{A}$-modules on $X$.
\end{lem}

By virtue of Proposition 11.1, see Mallios[\cite{mallios}, p.51],
if $\mathcal{E}_0$ and $\mathcal{F}_0$ are
sub-$\mathcal{A}$-modules of $\mathcal{E}$ and $\mathcal{F}$,
respectively, where $\mathcal{E}$ and $\mathcal{F}$ form a pairing
$[\mathcal{F}, \mathcal{E}; \mathcal{A}]$, then for any open
$U\subseteq X$, \[\begin{array}{ll} \mathcal{E}_0^\perp(U)=
\mathcal{E}_0(U)^\perp & \mbox{and $\mathcal{F}_0^\perp(U)=
\mathcal{F}_0(U)^\perp$}\end{array}\]up to
$\mathcal{A}(U)$-isomorphisms. It is however trivial that
$\mathcal{E}_0\subseteq (\mathcal{E}_0^\perp)^\perp:=
\mathcal{E}_0^{\perp\perp}$. Similarly, $\mathcal{F}_0\subseteq
\mathcal{F}_0^{\perp\perp}$. Of special importance is the
sub-$\mathcal{A}$-module $\mathcal{E}^\perp$ of $\mathcal{F}$,
that is the sub-$\mathcal{A}$-module \textit{orthogonal} to the
$\mathcal{A}$-module $\mathcal{E}$. We shall call
$\mathcal{E}^\perp$ the \textbf{left kernel}
sub-$\mathcal{A}$-module of the pairing $[\mathcal{F},
\mathcal{E}; \mathcal{A}]$. Similarly, we call $\mathcal{F}^\perp$
the \textbf{right kernel} sub-$\mathcal{A}$-module of
$[\mathcal{F}, \mathcal{E}; \mathcal{A}]$. Other authors such as
Crumeyrolle[\cite{rolle}, p.2] would call $\mathcal{E}^\perp$
(resp. $\mathcal{F}^\perp$) the \textbf{conjugate} of
$\mathcal{E}$ (resp. $\mathcal{F}$) in $\mathcal{F}$ (resp.
$\mathcal{E}$).

Also very important is our attempt, we are concerned with now, of
obtaining the \textit{kernel} of an $\mathcal{A}$-morphism $\phi:
\mathcal{E}\longrightarrow \mathcal{F}$ of $\mathcal{A}$-modules
($\mathcal{E}$ and $\mathcal{F}$ are defined on a topological space
$X$) as the sheafification of some presheaf of
$\mathcal{A}(U)$-modules on $X$. The kernel of $\phi$, denoted here
\[\mathcal{K}er \phi\equiv (\mathcal{K}er \phi, \pi|_{\mathcal{K}er
\phi}, X),\]is a sub-$\mathcal{A}$-module of $\mathcal{E}\equiv
(\mathcal{E}, \pi, X)$, see Mallios[\cite{mallios}, p.108].

\begin{lem}
Let $\mathcal{E}$ and $\mathcal{F}$ be $\mathcal{A}$-modules on a
topological space $X$, $\Gamma(\mathcal{E})\equiv
(\mathcal{E}(U)\equiv \Gamma(U, \mathcal{E}), \pi^U_V)$ the
$($complete$)$ presheaf of sections of $\mathcal{E}$, and
$\phi\equiv (\phi_U)\in Hom_\mathcal{A}(\mathcal{E}, \mathcal{F})$.
Furthermore, for every open subset $U\subseteq X$, let \[\ker
\phi_U= \{s\in \mathcal{E}(U):\ \phi_U(s)= 0\in
\mathcal{F}(U)\}.\]Then, the diagram
\[\xymatrix{U\ar[r] & \ker \phi_U \ar[d]^{\lambda^U_V} \\ V\ar[u]\ar[r] & \ker
\phi_V},\]where $\lambda^U_V:= \pi^U_V|_{\ker \phi_U}$, yields a
complete presheaf of $\mathcal{A}(U)$-modules, denoted \[\ker
\phi:= ((\ker \phi)(U)\equiv \ker \phi_U, \lambda^U_V),\] which is
the same as the $($complete$)$ presheaf of sections of the kernel
$\mathcal{K}er \phi\equiv (\mathcal{K}er \phi, \pi|_{\mathcal{K}er
\phi}, X)$ of $\phi$.
\end{lem}

\begin{proof}
That $\ker \phi\equiv (\ker \phi_U, \lambda^U_V)$ is a presheaf of
$\mathcal{A}(U)$-modules is easy to see.

Now, let $s\in (\mathcal{K}er \phi)(U)\equiv \Gamma(U,
\mathcal{K}er \phi)$, where $U$ is an open subset of $X$.
Evidently, for $x\in U$, $s(x)\in (\mathcal{K}er \phi)_x=
\mathcal{K}er \phi_x$, see Mallios[\cite{mallios}, (2.11), p.108].
We further obtain that, for all $x\in U$, \[\phi^\ast(s)(x):=
\phi_x(s(x))= 0\in \mathcal{F}_x,\]see Mallios[\cite{mallios},
Proposition 2.1, p.11]. It thus follows that $\phi^\ast(s)= 0$, or
equivalently $\phi_U(s):= \phi\circ s= \phi^\ast(s)= 0$. Hence,
$s\in \ker \phi_U$, so that $(\mathcal{K}er \phi)(U)\subseteq \ker
\phi_U$.

Conversely, let $s\in \ker \phi_U\subset \mathcal{E}(U)\equiv
\Gamma(U, \mathcal{E})$. Then, $\phi_x(s(x)):= \phi_U(s)(x)= 0$
for $x\in U$; consequently $s(x)\in (\mathcal{K}er \phi)_x$ for
$x\in U$. Since $s\in \mathcal{E}(U)$, it follows that $s\in
(\mathcal{K}er \phi)(U)$. Thus, $\ker \phi_U\subseteq
(\mathcal{K}er \phi)(U)$ for every open $U\subseteq X$; hence
$\ker \phi_U= (\mathcal{K}er \phi)(U)$ for every open $U\subseteq
X$.

We deduce from the above that $\ker \phi= \Gamma(\mathcal{K}er
\phi)$, and the proof is finished.
\end{proof}

We are ready now for one important result regarding pairings of
$\mathcal{A}$-modules. Let $\mathcal{E}$ and $\mathcal{F}$ be
$\mathcal{A}$-modules on a topological space $X$, and
\[\begin{array}{ll}\Gamma(\mathcal{E})\equiv (\Gamma(U, \mathcal{E}),
\pi^U_V) & \mbox{and $ \Gamma(\mathcal{F})\equiv (\Gamma(U,
\mathcal{F}), \rho^U_V)$}\end{array}\] their corresponding
(complete) presheaves of sections, respectively. We also assume
that $\Gamma(\mathcal{A})\equiv (\Gamma(U, \mathcal{A}),
\kappa^U_V)$ is the presheaf of sections for the sheaf
$\mathcal{A}$. Suppose that in the pairing $[\mathcal{F},
\mathcal{E}; \mathcal{A}]$ the \textit{left kernel} is $0$, i.e.
$\mathcal{E}^\perp =0$. Let $U$ be an open subset of $X$. For
every $r\in \mathcal{F}(U)$, consider the
$\mathcal{A}(U)$-morphism $\phi_r\equiv (\phi_{r, V})_{U\supseteq
V, \ open}\in \mathcal{E}^\ast(U)$, given by
\[\phi_r(t)= \rho^U_V(r)t\]for all $t\in \mathcal{E}(V)$, where $V$
is an open subset of $U$. One sees that the $\mathcal{A}(U)$-map
\[\begin{array}{ll}\Phi_U: \mathcal{F}(U)\longrightarrow \mathcal{E}^\ast(U); &
r\longmapsto \phi_r \end{array}\] is an
$\mathcal{A}(U)$-homomorphism. In fact, for all $t\in
\mathcal{E}(V)$, where $V$ is open in $U$, and all $a\in
\mathcal{A}(U)$ and $r, s\in \mathcal{F}(U)$, we have
\[\phi_{r+s}(t)= \rho^U_V(r+s)t = (\rho^U_V(r)+ \rho^U_V(s))t=
\rho^U_V(r)t+ \rho^U_V(s)t= \phi_r(t)+ \phi_s(t)\]and
\[\phi_{ar}(t)= \rho^U_V(ar)t= \kappa^U_V(a)\rho^U_V(r)t=
\kappa^U_V(a)\phi_r(t)\equiv a\phi_r(t).\]However, \[\ker \Phi_U=
\{r\in \mathcal{F}(U):\ \phi_r(t)=0\ \mbox{for all $t\in
\mathcal{E}(V)$, where $V$ is open in $U$}\}.\]Since we assumed
that $\mathcal{E}^\perp =0$, i.e. $\mathcal{E}^\perp(U)=
\mathcal{E}(U)^\perp= 0$ for all open $U\subseteq X$, it follows
that $r= 0$. Hence, for every open $U\subseteq X$, $\Phi_U$ is an
$\mathcal{A}(U)$-\textit{isomorphism into}. If we let $U$ vary
over the open subsets of $X$, the family $\Phi\equiv
(\Phi_U)_{X\supseteq U, \ open}$ is a
$\Gamma(\mathcal{A})$-morphism of presheaves $\Gamma(\mathcal{F})$
and $\Gamma(\mathcal{E}^\ast)$. In fact, first observe that if
$\Gamma(\mathcal{E}^\ast)\equiv (\mathcal{E}^\ast(U),
{\pi^\ast}^U_V)$, the restriction maps are defined as follows: For
$\alpha\equiv (\alpha_O)_{U\supseteq O,\ open}\in
\mathcal{E}^\ast(U):= Hom_{\mathcal{A}|_U}(\mathcal{E}|_U,
\mathcal{A}|_U),$ \[ {\pi^\ast}^U_V(\alpha):=
(\alpha_O)_{V\supseteq O,\ open}\in \mathcal{E}^\ast(V).\] Hence,
for $r\in \mathcal{F}(U)$ and $t\in (\mathcal{E}|_V)(W)=
\mathcal{E}(W)$, where $W$ is open in $V$, we have
\[\begin{array}{lll}\Phi_V(\rho^U_V(r))(t) & = &
\phi_{\rho^U_V(r)}(t) \\ & = & \rho^V_W(\rho^U_V(r))(t)\\ & = &
\rho^U_W(r)(t)\end{array}\]and \[\begin{array}{lll}
{\pi^\ast}^U_V(\Phi_U(r))(t) & = & {\pi^\ast}^U_V(\phi_r)(t)\\ & =
& (\phi_{r, O})_{V\supseteq O,\ open}(t)\\ & = &
\rho^U_W(r)(t);\end{array}\] thus the diagram
\[\xymatrix{\mathcal{F}(U)\ar[r]^{\Phi_U}\ar[d]_{\rho^U_V} & \mathcal{E}^\ast(U)\ar[d]^{{\pi^\ast}^U_V}\\
\mathcal{F}(V)\ar[r]_{\Phi_V} & \mathcal{E}^\ast(V)}\] commutes,
and consequently $\Phi\in
Hom_{\Gamma(\mathcal{A})}(\Gamma(\mathcal{F}),
\Gamma(\mathcal{E}^\ast)).$ Through the sheafification functor
$\mathbf{S}$, and on the basis of Mallios[\cite{mallios}, (13.19),
p.75], we have that
\[\mathbf{S}(\Phi)\equiv \widetilde{\Phi}\in
Hom_{\mathcal{S}h_X}(\mathcal{F}, \mathcal{E}^\ast)\]is an
$\mathcal{A}$-isomorphism into.

Similarly, suppose that $\mathcal{F}^\perp= 0$; let $\Psi_U$ be the
$\mathcal{A}(U)$-map, given by \[\begin{array}{ll}\Psi_U:
\mathcal{E}(U)\longrightarrow \mathcal{F}^\ast(U); & s\longmapsto
\Psi_U(s)\equiv \psi_s\end{array}\]where \[\psi_s(t)=
t\pi^U_V(s)\]for all $t\in \mathcal{F}(V)$ and an open $V\subseteq
U$. For $r, s\in \mathcal{E}(U)$, $a\in \mathcal{A}(U)$ and $t\in
\mathcal{F}(V)$, where as above $V$ is open in $U$, we have
\[\psi_{r+s}(t)= t(\pi^U_V(r)+ \pi^U_V(s))= t\pi^U_V(r)+
t\pi^U_V(s)= \psi_r(t)+ \psi_s(t)\]and \[\psi_{ra}(t)= t\pi^U_V(ra)=
t\pi^U_V(r)\kappa^U_V(a)= \psi_r(t)\kappa^U_V(a)\equiv
\psi_r(t)a;\]so that $\Psi_U$ is an $\mathcal{A}(U)$-homorphism.
Now, since $\mathcal{F}^\perp(U)= \mathcal{F}(U)^\perp= 0$ for every
open $U\subseteq X$,
\begin{eqnarray*}\ker \Psi_U:= \{s\in \mathcal{E}(U):\ \psi_s(t)= 0\ \mbox{for all $t\in \mathcal{F}(V)$
 and $V$ open in $U$}\}=0,\end{eqnarray*}for any open
$U\subseteq X$. Hence, every $\Psi_U$ is an
$\mathcal{A}(U)$-isomorphism into. As in the previous case, the
family $\Psi\equiv (\Psi_U)_{X\supseteq U,\ open}$ is a
$\Gamma(\mathcal{A})$-isomorphism into of presheaves
$\Gamma(\mathcal{E})$ and $\Gamma(\mathcal{F}^\ast)$.
Consequently,
\[\mathbf{S}(\Psi)\equiv \widetilde{\Psi}\in
Hom_{\mathcal{S}h_X}(\mathcal{E}, \mathcal{F}^\ast)\]is an
$\mathcal{A}$-isomorphism into.

Let us make a short breach here for the following useful lemmas.

\begin{lem}\label{lem1.3}
Let $(X, \mathcal{A})$ be a $\mathbb{C}$-algebraized space, and
$E\equiv (E(U), \rho^U_V)$ and $F\equiv (F(U), \sigma^U_V)$ be
presheaves of $\mathcal{A}(U)$-modules on $X$. Then,
\[\mathbf{S}(E\times F)= \mathbf{S}(E)\times \mathbf{S}(F)\]within
an $\mathcal{A}$-isomorphism.
\end{lem}

\begin{proof}
Sheaves $\mathbf{S}(E\times F)$ and $\mathbf{S}(E)\times
\mathbf{S}(F)$ clearly have the same stalks at every $x\in X$.
Therefore the underlying sets of $\mathbf{S}(E\times F)$ and
$\mathbf{S}(E)\times \mathbf{S}(F)$ are the same. It remains only
to show that the topology making $\mathbf{S}(E\times F)$ into a
sheaf is the same as the topology which defines the sheaf
structure on $\mathbf{S}(E)\times \mathbf{S}(F)$. The topology of
$\mathbf{S}(E\times F)$ is the topology generated by the basis
\[\{\tilde{s}(V):\ s\in (E\times F)(U),\
\mbox{where $U, V$ are open in $X$ with $V\subseteq U$}\},\] see
Mallios[\cite{mallios}, Theorem 3.1, p.14]. But $s\in (E\times
F)(U)= E(U)\times F(U)$ is if and only if it is of the form
\[s= (s_1, s_2)\]where $s_1\in E(U)$ and $s_2\in F(U)$. It follows
that \[\tilde{s}= (\widetilde{s_1},
\widetilde{s_2});\]consequently \[\tilde{s}(V)=
(\widetilde{s_1}(V), \widetilde{s_2}(V))\]for any open subset
$V\subseteq U\equiv Dom(\tilde{s})$. Besides,
\[\{\widetilde{s_1}(V):\ s_1\in E(U)\  \mbox{and $V$ is open in
$U$}\}\]and \[\{\widetilde{s_2}(V):\ s_1\in F(U)\  \mbox{and $V$ is
open in $U$}\}\] are bases for the topologies of $\mathbf{S}(E)$ and
$\mathbf{S}(F)$, respectively, therefore the topology of
$\mathbf{S}(E\times F)$ is equivalent to the topology of
$\mathbf{S}(E)\times \mathbf{S}(F)$; thus the proof is finished.
\end{proof}

\begin{lem}\label{lem1.4}
Let $E\equiv (E(U), \rho^U_V)$, $F\equiv (F(U), \sigma^U_V)$ and
$A\equiv (A(U), \kappa^U_V)$ be presheaves of $A(U)$-modules on a
topological space $X$. Suppose that a map $\phi\in
Hom_{\mathcal{P}\mathcal{S}h_X}(E\times F, A)$ is given, and
$\mathcal{E}:= \mathbf{S}(E)$, $\mathcal{F}:= \mathbf{S}(F)$, and
$\mathcal{A}:= \mathbf{S}(A)$, where $\mathbf{S}:
\mathcal{P}\mathcal{S}h_X\longrightarrow \mathcal{S}h_X$ is the
sheafification functor. Then, if $\phi$ is bilinear, i.e. every
$\phi_U: E(U)\times F(U)\longrightarrow A(U)$, where $U$ is open in
$X$, is bilinear,  the $\mathcal{A}$-morphism
$\mathbf{S}(\phi)\equiv \overline{\phi}\in
Hom_{\mathcal{S}h_X}(\mathcal{E}\times \mathcal{F}, \mathcal{A})$ is
also bilinear, and
\begin{equation}\label{eq4}\overline{\phi}_U(\tilde{s}, \tilde{t})=
\widetilde{\phi_U(s, t)}\end{equation}for $s\in E(U)$, $t\in F(U)$,
$\tilde{s}=\rho_U(s)$, and $\tilde{t}= \sigma_U(t)$, cf.
Mallios$[\cite{mallios}, (7.22), p. 32].$
\end{lem}

\begin{proof}
First, let us make the following comment. Given a presheaf $S\equiv
(S(U), \rho^U_V)$ on a topological space $X$, the sheafification of
$S$ hinges for every open subset $U$ of $X$ on the corresponding map
$\rho_U: S(U)\longrightarrow \Gamma(U, \mathcal{S})\equiv
\mathcal{S}(U)$, which associates with every $s\in S(U)$ the section
$\rho_U(s)\equiv\tilde{s}\in \mathcal{S}(U)$. In the special case
where $S$ is a presheaf of $A(U)$-modules, maps $\rho_U$ are
$A(U)$-homomorphisms, so that $\rho(s+t)\equiv \widetilde{s+t}=
\tilde{s}+ \tilde{t}\equiv \rho_U(s)+ \rho_U(t)$, for all $s, t\in
S(U)$. Thus, we have the following.

Given that $\mathbf{S}$ is a functor, the diagram
\[\xymatrix{E(U)\times F(U)\ar[r]^{\rho_U\times \sigma_U}\ar[d]_{\phi_U} & \mathcal{E}(U)\times
\mathcal{F}(U)\ar[d]^{\mathbf{S}(\phi)\equiv \overline{\phi}} \\
A(U)\ar[r]_{\kappa_U} & \mathcal{A}(U)},\] where $\rho_U:
E(U)\longrightarrow \mathcal{E}(U)$, $\sigma_U: F(U)\longrightarrow
\mathcal{F}(U)$, and $\kappa_U: A(U)\longrightarrow \mathcal{A}(U)$
are the (canonical) maps defining the respective sheafifications, is
commutative, and one has Equation (\ref{eq4}). For $s, s'\in E(U)$
and $t\in F(U)$, where $U$ is an open subset of $X$, it is easy to
see that
\[\mathbf{S}(\phi_U)(\widetilde{s+ s'}, \tilde{t})=
\mathbf{S}(\phi_U)(\tilde{s}, \tilde{t})+
\mathbf{S}(\phi_U)(\widetilde{s'}, \tilde{t}).\] Likewise one shows
linearity in the second component, and that is the end of the proof.
\end{proof}

Now, let us return to the assumption $\mathcal{E}^\perp= 0$, and
let $\mathcal{F}_0$ be a sub-$\mathcal{A}$-module of
$\mathcal{F}$. Following Artin[\cite{artin}, p.19], we find in a
natural way a new pairing: more precisely, the \textit{pairing of
$\mathcal{A}$-modules $\mathcal{F}_0$} and
\textit{$\mathcal{E}/\mathcal{F}_0^\perp$ into $\mathcal{A}$.} For
this purpose, let $U$ be an open subset of $X$; given $t\in
\mathcal{F}_0(U)$ and $s\in \mathcal{E}(U)$, we define as product
of $t$ and $s+ \mathcal{F}_0^\perp(U)$ the element $ts$ of
$\mathcal{A}(U)$: \begin{equation}\label{eq3}t\cdot(s+
\mathcal{F}_0^\perp(U)):= ts.
\end{equation}That this multiplication is well defined is easy to
see. In fact, suppose $s+ \mathcal{F}_0^\perp(U)= s_1+
\mathcal{F}_0^\perp(U)$; therefore $s-s_1\in
\mathcal{F}_0^\perp(U)$. But $t\in \mathcal{F}_0(U)$, so
$t(s-s_1)=0$, i.e. $ts= ts_1$. It is also easy to see that the
multiplication (\ref{eq3}) is bilinear, thus we have the pairing
$[\mathcal{F}_0(U), \mathcal{E}(U)/\mathcal{F}_0^\perp(U);
\mathcal{A}(U)]$. In turn, if we let $U$ run over the open subsets
of $X$ and every multiplication $\mathcal{F}_0(U)\times
\mathcal{E}(U)/\mathcal{F}_0^\perp(U)\longrightarrow \mathcal{A}(U)$
is given as in Equation (\ref{eq3}), we obtain a
$\Gamma(\mathcal{A})$-morphism \[\xymatrix{\mathcal{F}_0(U)\times
(\mathcal{E}(U)/\mathcal{F}_0^\perp(U))\ar[r]\ar[d] &
\mathcal{A}(U)\ar[d]
\\ \mathcal{F}_0(V)\times
(\mathcal{E}(V)/\mathcal{F}_0^\perp(V))\ar[r] &
\mathcal{A}(V)}\]of $\Gamma(\mathcal{A})$-presheaves
$\Gamma(\mathcal{F}_0)\times
\Gamma(\mathcal{E})/\Gamma(\mathcal{F}_0^\perp)$ and
$\Gamma(\mathcal{A})$, see Mallios[\cite{mallios}, pp 114-115].
But, by virtue of Lemma \ref{lem1.3},
$$\mathbf{S}(\Gamma(\mathcal{F}_0)\times
\Gamma(\mathcal{E})/\Gamma(\mathcal{F}_0^\perp))=
\mathcal{F}_0\times (\mathcal{E}/\mathcal{F}_0^\perp)$$ within an
$\mathcal{A}$-isomorphism, it follows through the sheafification
functor and from Lemma \ref{lem1.4} that the
$\mathcal{A}$-morphism
$$\mathcal{F}_0\times
(\mathcal{E}/\mathcal{F}_0^\perp)\longrightarrow \mathcal{A}$$ is
bilinear. Hence, we have a pairing $[\mathcal{F}_0,
\mathcal{E}/\mathcal{F}_0^\perp; \mathcal{A}]$. Now, what is the
\textit{right} kernel of the pairing $[\mathcal{F}_0,
\mathcal{E}/\mathcal{F}_0^\perp; \mathcal{A}]$? We shall denote
this right kernel by \[\widehat{\mathcal{F}_0}\]to differentiate
it from the orthogonal $\mathcal{F}_0^\perp$. First, let us
observe the following. Let $U$ be an open subset of $X$; the right
kernel of the pairing $[\mathcal{F}_0(U),
\mathcal{E}(U)/\mathcal{F}_0^\perp(U); \mathcal{A}(U)]$, which we
denote by $\widehat{\mathcal{F}_0(U)}$, consists of elements $s+
\mathcal{F}_0^\perp(U)$ such that $t\cdot(s+
\mathcal{F}_0^\perp(U)):= ts=0$ for all $t\in \mathcal{F}_0(U)$.
This means that $s\in \mathcal{F}_0(U)^\perp=
\mathcal{F}_0^\perp(U)$ and, therefore $s+ \mathcal{F}_0^\perp(U)=
\mathcal{F}_0^\perp(U)$. Thus, $\widehat{\mathcal{F}_0(U)}= 0$.

Now, we will show that $\widehat{\mathcal{F}_0}=0$. Indeed, let
$\tilde{s}\in \widehat{\mathcal{F}_0}(U)\subseteq
(\mathcal{E}/\mathcal{F}_0^\perp)(U)$. There exists $s+
\mathcal{F}_0^\perp(U)$ such that (see Mallios[\cite{mallios},
(7.9), p.30]) \[\tilde{s}(x)= [s+ \mathcal{F}_0^\perp(U)]_x\]for
$x\in U$. Let $t\in \mathcal{F}_0(U)$. We have \[(t\tilde{s})(x)=
t(x)\tilde{s}(x)= t(x)[s+ \mathcal{F}_0^\perp(U)]_x= 0;\]since $t$
is arbitrary, $s\in \mathcal{F}_0^\perp(U)$, so that $\tilde{s}= 0$.
Therefore, the right kernel of our pairing is $0$, as desired, and
$\widehat{\mathcal{F}_0}$ is the sheafification of the presheaf
\[\xymatrix{U\ar[r] & \widehat{\mathcal{F}_0(U)}\ar[d] \\ V\ar[u]\ar[r] & \widehat{\mathcal{F}_0(V)}}.\]

Furthermore, for every open $U\subseteq X$, we can construct,
using the previously established method, an
$\mathcal{A}(U)$-isomorphism into
\[\mathcal{E}(U)/\mathcal{F}_0^\perp(U)\longrightarrow
\mathcal{F}_0^\ast(U);\]consequently through the sheafification
functor we have an $\mathcal{A}$-isomorphism into:
\[\mathcal{E}/\mathcal{F}_0^\perp\longrightarrow
\mathcal{F}_0^\ast.\]

If $\mathcal{E}_0$ is a given sub-$\mathcal{A}$-module of
$\mathcal{E}$, we can also define a natural pairing
$[\mathcal{E}_0^\perp, \mathcal{E}/\mathcal{E}_0; \mathcal{A}]$ by
setting \[\begin{array}{lll} t\cdot(s+ \mathcal{E}_0(U))=ts, &
t\in \mathcal{E}_0^\perp(U)= \mathcal{E}_0(U)^\perp, & s+
\mathcal{E}_0(U)\in
\mathcal{E}(U)/\mathcal{E}_0(U).\end{array}\]Likewise, we obtain
an $\mathcal{A}$-isomorphism into:
\[\mathcal{E}_0^\perp\longrightarrow
(\mathcal{E}/\mathcal{E}_0)^\ast.\]

We may formulate our results as follows.

\begin{theo}\label{theo2}
Let $\mathcal{F}$ and $\mathcal{E}$ be $\mathcal{A}$-modules on a
topological space $X$ paired into a $\mathbb{C}$-algebra sheaf
$\mathcal{A}$, and assume that $\mathcal{E}^\perp=0$. Moreover,
let $\mathcal{F}_0$ be a sub-$\mathcal{A}$-module of $\mathcal{F}$
and $\mathcal{E}_0$ a sub-$\mathcal{A}$-module of $\mathcal{E}$.
There exist natural $\mathcal{A}$-isomorphisms into:
\[\begin{array}{lll}\mathcal{E}/\mathcal{F}_0^\perp\longrightarrow
\mathcal{F}_0^\ast, &  \mbox{and} &
\mathcal{E}_0^\perp\longrightarrow (\mathcal{E}/\mathcal{E}_0)^\ast.
\end{array}\]
\end{theo}

Now, as a special case assume in the pairing $[\mathcal{F},
\mathcal{E}; \mathcal{A}]$ of Theorem \ref{theo2} that both
$\mathcal{F}$ and $\mathcal{E}$ are free $\mathcal{A}$-modules of
finite rank $m$ and $n$ respectively. Suppose that $\phi:
\mathcal{F}\times \mathcal{E}\longrightarrow \mathcal{A}$ is the
bilinear morphism which defines the pairing. Define
$\mathcal{A}$-morphisms $\gamma: \mathcal{E}\longrightarrow
\mathcal{F}^\ast= \mathcal{F}$ and $\delta:
\mathcal{F}\longrightarrow \mathcal{E}^\ast= \mathcal{E}$ (see
Mallios[\cite{mallios}, (5.2), p. 298]) such that
\[\gamma_U(s)(t)\equiv \gamma^s_U(t):= \phi_U(t, s)\]and
\[\delta_U\equiv \delta^t_U(s):= \phi_U(t, s)\]for all $s\in
\mathcal{E}(U)$ and $t\in \mathcal{F}(U)$, where $U$ is an open
subset of $X$. For every open subset $U$ of $X$, $\ker\gamma_U=
\mathcal{F}^\perp(U)$ and $\ker\delta_U= \mathcal{E}^\perp(U)$;
therefore \[\begin{array}{ll} \mathcal{K}er\gamma=
\mathcal{F}^\perp, & \mathcal{K}er\delta=
\mathcal{E}^\perp\end{array}\]within $\mathcal{A}$-isomorphisms
respectively. The $\mathcal{A}$-morphism $\phi$ is said to be
\textbf{non-degenerate} if $\mathcal{E}^\perp= \mathcal{F}^\perp=
0$, and \textbf{degenerate} otherwise.

If $\mathcal{E}^\perp= \mathcal{F}^\perp= 0$, $\gamma$ and
$\delta$ are injective, or equivalently $\gamma_x:
\mathcal{E}_x\longrightarrow \mathcal{F}_x$ and $\delta_x:
\mathcal{F}_x\longrightarrow \mathcal{E}_x$ are injective for
every $x\in X$, cf. Mallios[\cite{mallios}, Lemma 12.1, p.60]. It
follows clearly that $m= n$, i.e. $\mathcal{F}$ is
$\mathcal{A}$-isomorphic to $\mathcal{E}$.

On the other hand, suppose that $\mathcal{E}^\perp$ and
$\mathcal{F}^\perp$ are not all zero. For this case, we will need
the following lemma.

\begin{lem}\label{lem1.5}
Let $A\equiv (A(U), \kappa^U_V)$ be a presheaf of
$\mathbb{C}$-algebras on a topological space $X$, $E\equiv (E(U),
\rho^U_V)$ and $F\equiv (F(U), \sigma^U_V)$ $A$-presheaves $($i.e.
presheaves of $A(U)$-modules$)$ on $X$, and finally $\Phi\equiv
(\Phi_U)_{X\supseteq U,\ open}$ a bilinear $A$-morphism $\Phi:
F\times E\longrightarrow A$. For every open $U\subseteq X$, let
\[E(U)^\perp\equiv E^\perp(U):= \{t\in F(U):\ \Phi_U(t, s)= 0\
\mbox{for all $s\in E(U)$}\}.\]Assume that $E^\perp\equiv
(E^\perp\equiv E(U)^\perp, \lambda^U_V\equiv
\sigma^U_V|_{E^\perp(U)})$ is a presheaf on $X$. Then, if
$\mathcal{E}\equiv \mathbf{S}(E)$ and $\mathcal{A}\equiv
\mathbf{S}(A)$, one has that
\[\mathcal{E}^\perp= \mathbf{S}(E^\perp)\]within an
$\mathcal{A}$-isomorphism.
\end{lem}

\begin{proof}
Let $U$ be an open subset of $X$; $\mathbf{S}(E^\perp)(U)$
consists of elements (in fact (local) sections) $\tilde{t}\equiv
\rho_U(t)$, where $t\in E^\perp(U)$ and $\rho_U:
E^\perp(U)\longrightarrow \Gamma(U, \mathbf{S}(E^\perp))\equiv
\mathbf{S}(E^\perp)(U)$ is the canonical map obtained through the
sheafification process, see Mallios[\cite{mallios}, (7.22), p.
32]. Let $\mathcal{F}\equiv \mathbf{S}(F)$, and denote by
$\overline{\Phi}$ the bilinear $\mathcal{A}$-morphism
$\overline{\Phi}: \mathcal{F}\times \mathcal{E}\longrightarrow
\mathcal{A}$ induced by the presheaf morphism $\Phi$, see Lemma
\ref{lem1.4}. But, for all $t\in E^\perp(U)$ and $s\in E(U)$, we
have \[\overline{\Phi}_U(\tilde{s}, \tilde{t})=
\widetilde{\Phi_U(s, t)}= 0,\]which means that $\tilde{t}\in
\mathbf{S}(E)(U)$ if and only if $\tilde{t}\in
\mathcal{E}^\perp(U)$; hence $\mathbf{S}(E)(U)=
\mathcal{E}^\perp(U)$, and the proof is finished.
\end{proof}

Getting back to the case where $\mathcal{E}^\perp$ and
$\mathcal{F}^\perp$ are not all zero in the pairing $[\mathcal{F},
\mathcal{E}; \mathcal{A}]$, defined by the map $\Phi:
\mathcal{F}\times \mathcal{E}\longrightarrow \mathcal{A}$, and
where $\mathcal{F}$ and $\mathcal{E}$ are free
$\mathcal{A}$-modules of finite rank on $X$, we notice that for
every open subset $U$ of $X$, if $t, t_1\in \mathcal{F}(U)$, and
$t-t_1\in \mathcal{E}^\perp(U)$, then $\Phi_U(t, s)= \Phi(t_1, s)$
for all $s\in \mathcal{E}(U)$. Analogously if $s, s_1\in
\mathcal{E}(U)$ and $t, t_1\in \mathcal{F}(U)$ such that $s-s_1\in
\mathcal{F}^\perp(U)$ and $t-t_1\in \mathcal{E}^\perp(U)$, we have
\[\Phi_U(t, s)= \Phi_U(t_1, s_1).\]Setting \[\begin{array}{ll}\bar{s}\equiv
[s] (\mbox{mod}\ \mathcal{F}^\perp(U)), & \bar{t}= [t]
(\mbox{mod}\ \mathcal{E}^\perp(U))\end{array}\]we obtain the
bilinear $\mathcal{A}(U)$-morphism \[\overline{\Phi}_U:
\mathcal{F}(U)/\mathcal{E}^\perp(U)\times
\mathcal{E}(U)/\mathcal{F}^\perp(U)\longrightarrow
\mathcal{A}(U),\]given by \[\overline{\Phi}_U(\bar{t}, \bar{s})=
\Phi_U(t, s)\]for all $t\in \mathcal{F}(U)$ and $s\in
\mathcal{E}(U)$. Clearly, that $\overline{\Phi}_U(\bar{t},
\bar{s})= 0$ for all $\bar{s}\in
\mathcal{E}(U)/\mathcal{F}^\perp(U)$ is equivalent to saying that
$\Phi_U(t, s)= 0$ for all $s\in \mathcal{E}(U)$, so that $t\in
\mathcal{E}^\perp(U)$, and $\bar{t}=0$. Similarly,
$\overline{\Phi}_U(\bar{t}, \bar{s})= 0$ for all $\bar{t}\in
\mathcal{F}(U)/\mathcal{E}^\perp(U)$ is equivalent to $\bar{s}=0$.
Hence, $(\mathcal{F}(U)/\mathcal{E}^\perp(U))^\perp= 0$ and
$(\mathcal{E}(U)/\mathcal{F}^\perp(U))^\perp=0$, i.e.
$\overline{\Phi}_U$ is non-degenerate. By Lemma \ref{lem1.5},
\[(\mathcal{F}/\mathcal{E}^\perp)^\perp\equiv
\mathbf{S}(\Gamma(\mathcal{F})/\Gamma(\mathcal{E}^\perp))^\perp=
\mathbf{S}((\Gamma(\mathcal{F})/\Gamma(\mathcal{E}^\perp))^\perp);\]since
$(\mathcal{F}(U)/\mathcal{E}^\perp(U))^\perp=0$ for all open
subset $U\subseteq X$, it follows that
$(\mathcal{F}/\mathcal{E}^\perp)^\perp=0.$ Likewise
$(\mathcal{E}/\mathcal{F}^\perp)^\perp=0$. Thus, we have

\begin{theo}\label{theo3}
Let $(X, \mathcal{A})$ be a $\mathbb{C}$-algebraized space,
$\mathcal{F}$ and $\mathcal{E}$ be free $\mathcal{A}$-modules of
finite rank, paired into $\mathcal{A}$ through a bilinear morphism
$\Phi: \mathcal{F}\times \mathcal{E}\longrightarrow \mathcal{A}$.
The following hold: \begin{enumerate}\item[{$(i)$}] If $\Phi$ is
non-degenerate, then $\mathcal{F}= \mathcal{E}$ within an
$\mathcal{A}$-isomorphism. \item[{$(ii)$}] If $\Phi$ is
degenerate, then $\mathcal{F}/\mathcal{E}^\perp$ and
$\mathcal{E}/\mathcal{F}^\perp$ are paired into $\mathcal{A}$
through the $\mathcal{A}$-morphism $\overline{\Phi}:
\mathcal{F}/\mathcal{E}^\perp\times
\mathcal{E}/\mathcal{F}^\perp\longrightarrow \mathcal{A},$ given
by \[\overline{\Phi}_U(\bar{t}, \bar{s})= \Phi_U(t, s)\]for all
$t\in \mathcal{F}(U)$, $s\in \mathcal{E}(U)$, and where $U$ is an
arbitrary open subset of $X$. Moreover, $\overline{\Phi}$ is
non-degenerate, i.e. \[\begin{array}{ll}
(\mathcal{F}/\mathcal{E}^\perp)^\perp=0, &
(\mathcal{E}/\mathcal{F}^\perp)^\perp=0.
\end{array}\]\end{enumerate}
\end{theo}

\noindent Anastasios Mallios\\ {Department of Mathematics} \\
{University of Athens}\\ {Athens, Greece}\\
{Email: amallios@math.uoa.gr}

 \noindent Patrice P.
Ntumba\\{Department of Mathematics and Applied
Mathematics}\\{University of Pretoria}\\ {Hatfield 0002, Republic
of South Africa}\\{Email: patrice.ntumba@up.ac.za}

\end{document}